\documentclass[12pt]{article}
\usepackage[intlimits]{amsmath}
\usepackage{amsfonts,amssymb,amscd,amsthm,cite,hyperref}
\input xy
\xyoption{all}

\usepackage{ifpdf}
\ifpdf 
  \usepackage[pdftex]{graphicx}
  \DeclareGraphicsExtensions{.pdf,.png,.jpg,.jpeg,.mps}
  \usepackage{pgf}
  \usepackage{tikz}
\else 
  \usepackage{graphicx}
  \DeclareGraphicsExtensions{.eps,.bmp}
  \usepackage{pgf}
  \usepackage{tikz}
  \usepackage{pstricks}
\fi

\def\ind{\operatorname{ind}}

\def\Td{\operatorname{Td}}
\def\ch{\operatorname{ch}}

\def\Mat{\operatorname{Mat}}
\DeclareMathOperator\tr{Tr}

\def\dim{\operatorname{dim}}

\def\det{\operatorname{det}\nolimits}

\newtheorem{proposition}{Proposition}
\newtheorem{theorem}{Theorem}

\newtheorem{lemma}{Lemma}

\theoremstyle{definition}

\newtheorem{definition}{Definition}

\title{Index of nonlocal elliptic boundary value problems associated with isometric group actions}

\author{Anton Savin}

\date{}
\begin{document}

\maketitle

\begin{abstract} 
Given a compact manifold   with boundary endowed with an isometric action of a discrete group of polynomial growth, we state an index theorem for elliptic elements in the algebra of nonlocal operators generated by the Boutet de Monvel algebra of pseudodifferential boundary value problems on the manifold and the shift operators associated with the group action.
\end{abstract}



\section{Introduction}

Let $\Gamma$ be a discrete group of diffeomorphisms of a smooth manifold $M$. We consider the class of operators generated by pseudodifferential operators on $M$ and shift operators $T_\gamma u(x)=u(\gamma^{-1}(x))$ for all $\gamma\in\Gamma$, while $u\in C^\infty(M)$. Fredholm property for operators in this class is known in quite general situation (see \cite{AnLe1,AnLe2}). However,  the index problem  was only studied in the case of manifolds without boundary (see \cite{Ant4,NaSaSt17,SaSt39,SaSt40,Per5} and the references cited there). 

In this paper, we consider  a compact smooth manifold   with boundary endowed with an isometric action of a discrete group of polynomial growth in the sense of Gromov \cite{Gro1}. In this geometric situation, we state an index theorem for elliptic elements in the algebra generated by the Boutet de Monvel    pseudodifferential boundary value problems \cite{Bout2} on the manifold and the shift operators associated with the group action. Our index formula gives as special cases the   index formula for elliptic elements in the Boutet de Monvel algebra \cite{Fds2} (see also \cite{MeScSc1}) and the index formula for elliptic operators associated with isometric group actions on   closed manifolds~\cite{NaSaSt17}.

The work was partially supported by RFBR grant~No~19-01-00447 and RUDN University program 5-100.

\section{Ellipticity and Fredholm property}

\paragraph{Boutet de Monvel algebra.}

Let  $M$ be a compact smooth manifold with boundary denoted by $X$. 
We suppose that $M$ is endowed with a Riemannian metric and consider the induced Riemannian metric on  $X$. 
Local coordinates on $M$ and $X$ are denoted by $x$ and $x'$, respectively. In addition,  in a neighborhood of the boundary we use coordinates   $x=(x',x_n),$ $x_n\ge 0$,
such that the boundary has the equation $x_n=0$, while $x_n$ is equal to the distance to the boundary.
The dual coordinates in $T^*M$ are denoted by $\xi=(\xi',\xi_n).$

We consider Boutet de Monvel operators of order and type equal to zero. We refer the reader to \cite{Bout2,ReSc1,Gru3,Schr3} for a complete exposition of the Boutet de Monvel algebra and recall here only several facts about this algebra, which are used below. 

Boutet de Monvel operators of order and type equal to zero define continuous mappings of the form
\begin{equation}
\label{eq-opb1}
\mathcal{D}= \left(
 \begin{array}{cc}
   A+G & C\\ 
   B & A_X 
 \end{array}
\right):
\begin{array}{c}
   L^2(M)\\ 
   \oplus \\
   L^2(X)
 \end{array}
 \longrightarrow
 \begin{array}{c}
   L^2(M)\\ 
   \oplus \\
   L^2(X)
 \end{array}
\end{equation}
Here
\begin{itemize}
\item $A$ is a pseudodifferential operator of order zero on $M$, whose complete symbol satisfies the so-called transmission property;
\item $A_X$ is a  pseudodifferential operator of order zero on $X$;
\item $B,C$ and $G$ are  boundary (or trace), coboundary (or potential)  and Green operators, respectively.
\end{itemize}

{\em The symbol} of   operator \eqref{eq-opb1} is a pair $\sigma(\mathcal{D})=(\sigma_M(\mathcal{D}),\sigma_X(\mathcal{D}))$. Here the first component is called the {\em interior symbol}  and is a function  
\begin{equation}\label{s-bdm0}
 \sigma_M(\mathcal{D})=\sigma(A)\in C^\infty(T^*_0M)
\end{equation}
homogeneous  on the cotangent bundle $T^*_0M$ with zero section deleted and   equal to the principal symbol of $A$. 
The second component is called  the {\em boundary symbol} and   is an operator function
$$
 \sigma_X(\mathcal{D}) \in C^\infty(T^*_0X,\mathcal{B}(\overline H_+\oplus \mathbb{C}))\simeq 
 C^\infty(T^*_0X,\mathcal{B}(L^2(\mathbb{R}_+)\oplus \mathbb{C})).
$$
Here $\overline H_+\subset L^2(\mathbb{R}_{\xi_n})$ is the Fourier image of the subspace $L^2(\mathbb{R}_+)\subset L^2(\mathbb{R}_{x_n})$ of functions vanishing for $x_n\le 0$.  The boundary symbol is twisted homogeneous in the following sense
\begin{equation}\label{eq-twist1}
\sigma_X(\mathcal{D})(x',\lambda \xi')=\varkappa^{-1}_\lambda \sigma_X(\mathcal{D})(x',  \xi') \varkappa_\lambda \quad
\text{for all }(x',\xi')\in T^*_0X,\lambda>0
\end{equation}
with respect to the unitary representation of $\mathbb{R}_+$ on $\overline H_+\oplus \mathbb{C}$:
$$
\varkappa_\lambda(u(\xi_n),v)=(\lambda^{1/2}u(\lambda\xi_n),v).
$$

We denote the algebra of boundary symbols by $\Sigma_X$. Let us describe elements in this algebra.
To this end, let $\theta_\pm(x_n)$ be the  functions on $\mathbb{R}$ equal to $1$ for $x_n\in \mathbb{R}_\pm$
and   zero otherwise. Consider Fr\'echet spaces
$
H_\pm=\mathcal{F}_{x_n\to\xi_n}(\theta_\pm(x_n)S(\mathbb{R})),
$ 
where $S(\mathbb{R})$ is the Schwarz space on $\mathbb{R}$. We also define
the projection $\Pi_+:H_+\oplus H_-\to H_+$ and the continuous functional
$$
\begin{array}{ccc}
\Pi':H_+\oplus H_- & \longrightarrow & \mathbb{C}\\
u(\xi_n) & \longmapsto &\lim\limits_{x_n\to 0+} \mathcal{F}^{-1}_{\xi_n\to x_n}(u(\xi_n)).
\end{array}
$$
Note that $\Pi'$ is just the Fourier image of the right limit at zero. Clearly, if $u\in L^1$, then
$$
\Pi'u=\frac{1}{2\pi}\int_\mathbb{R}u(\xi_n)d\xi_n.
$$
Let us now describe the elements in $\Sigma_X$. Consider smooth functions
\begin{itemize}
\item $b(x',\xi',\xi_n)\in C^\infty(T^*_0X,H_-)$;
\item $c(x',\xi',\xi_n)\in C^\infty(T^*_0X,H_+)$;
\item $g(x',\xi',\xi_n,\eta_n)\in C^\infty(T^*_0X,H_+\otimes H_-)$ (here we consider the topological tensor product of locally convex linear topological spaces $H_\pm$);
\item $q(x',\xi')\in C^\infty(T^*_0X)$.
\end{itemize}
We use these functions to define the family of operators
\begin{equation}
\label{eq-spb1}
a_X= \left(
 \begin{array}{cc}
   \Pi_+a(x',0,\xi',\xi_n)+\Pi'_{\eta_n}g(x',\xi',\xi_n,\eta_n) & c(x',\xi',\xi_n)\\ 
   \Pi'_{\xi_n}b(x',\xi',\xi_n) &  q(x',\xi') 
 \end{array}
\right):
\begin{array}{c}
   \overline{H}_+\\ 
   \oplus \\
   \mathbb{C}
 \end{array}
 \longrightarrow
 \begin{array}{c}
   \overline{H}_+\\ 
   \oplus \\
   \mathbb{C}
 \end{array}
\end{equation}
parametrized by $(x',\xi')\in T^*_0X$. Here $a(x',0,\xi',\xi_n)$ is the restriction of a zero-order symbol
on $M$ with the transmission property to the boundary and is called the {\em principal symbol} of   $a_X$. We suppose the functions $b,c,g,q$ in \eqref{eq-spb1} are chosen in such a way that the symbol $a_X$ is twisted-homogeneous (see~\eqref{eq-twist1}).
Then $a_X\in\Sigma_X$ and all elements in $\Sigma_X$ can be written as in \eqref{eq-spb1}. 

Finally, we define the linear mapping (regularized trace):
\begin{equation}\label{eq-regtr1}
\begin{array}{ccc}
\tr':\Sigma_X & \longrightarrow & C^\infty(T^*_0X)\\
\left(
 \begin{array}{cc}
   \Pi_+a +\Pi' g  & c \\ 
   \Pi' b  &  q 
 \end{array}
\right) & \longmapsto & \Pi'_{\xi_n}g(x',\xi',\xi_n,\xi_n)+q(x',\xi').
\end{array}
\end{equation}
This functional does not vanish on commutators in general as the following lemma shows.
\begin{lemma}[\cite{Fds2}]
Given $a_{1,X},a_{2,X}\in\Sigma_X$, we have
\begin{equation}\label{eq-regtr2}
 \tr'[a_{1,X},a_{2,X}]=-i\Pi'\left(\frac{\partial a_1}{\partial\xi_n}a_2\right)=
 i\Pi'\left(a_1\frac{\partial a_2}{\partial\xi_n}\right),
\end{equation}
where $a_{1,2}$ are the principal symbols of $a_{1,2,X}.$
\end{lemma}

\paragraph{$\Gamma$-Boutet de Monvel operators. Fredholm property.}

Let $\Gamma$ be a discrete group of isometries of   $M$.  Given $\gamma\in \Gamma$, we define shift operator
$$
T_\gamma:L^2(M) \oplus L^2(X)\longrightarrow L^2(M) \oplus L^2(X),\quad (u(x),v(x'))\longmapsto (u(\gamma^{-1}(x)),v(\gamma^{-1}(x'))).
$$
This operator is unitary if we endow  the $L^2$-spaces with the norms defined by the Riemannian volume forms associated with the metrics.  The mapping $\gamma\mapsto T_\gamma$ defines a unitary representation of $\Gamma$ on $L^2(M) \oplus L^2(X)$.

\begin{definition}
A $\Gamma$-{\em Boutet de Monvel operator}  is an operator equal to the sum
\begin{equation}\label{eq-op1}
 \mathcal{D}=\sum_{\gamma\in \Gamma} \mathcal{D}_\gamma T_\gamma:L^2(M) \oplus L^2(X)\longrightarrow L^2(M) \oplus L^2(X),
\end{equation}
where $\{\mathcal{D}_\gamma\}_{\gamma\in \Gamma}$ are Boutet de Monvel operators. We suppose that the sum in \eqref{eq-op1} is finite, i.e., only a finite number of $\mathcal{D}_\gamma$'s is nonzero.    
\end{definition}
Later on $\Gamma$-Boutet de Monvel operators are referred to as $\Gamma$-operators for short. 
One can show that given a Boutet de Monvel operator $\mathcal{D}$ and $\gamma\in \Gamma$ the composition
$T_\gamma\mathcal{D}T_\gamma^{-1}$ is also a Boutet de Monvel operator. This implies that operators  \eqref{eq-op1}
form an algebra. Moreover,  the interior and boundary symbols of $T_\gamma\mathcal{D}T_\gamma^{-1}$ are equal to
$$
\sigma_M(T_\gamma\mathcal{D}T_\gamma^{-1})(x,\xi)=\sigma_M(\mathcal{D})(\partial\gamma^{-1}(x,\xi)),
\quad \sigma_X(T_\gamma\mathcal{D}T_\gamma^{-1})(x',\xi')=\sigma_X(\mathcal{D})(\partial\gamma^{-1}(x',\xi')).
$$
Here the action of $\Gamma$ on $M$ and $X$ is lifted to the bundles $T^*M$ and $T^*X$  by the codifferentials
$ 
\partial \gamma=(d\gamma^t)^{-1} 
$
of the corresponding diffeomorphisms.  

We consider smooth crossed products (see \cite{Schwe1}) of algebras of interior and boundary symbols with  $\Gamma$ acting on these algebras by automorphisms. Recall that the smooth crossed product $\mathcal{A}\rtimes \Gamma$ of a Fr\'echet algebra $\mathcal{A}$ with the seminorms $\|\cdot\|_m$, $m> 0$, and a group $\Gamma$ of polynomial growth acting on $\mathcal{A}$ by automorphisms $a\mapsto \gamma(a)$ for all $a\in \mathcal{A}$ and $\gamma\in\Gamma$ is equal to the vector space of  functions $f:\Gamma\to \mathcal{A}$, which decay rapidly at infinity in the sense that the following estimates are valid:
$$
\|f(\gamma)\|_m\le C_N(1+|\gamma|)^{-N}
$$
for all $N,m>0$ and $\gamma\in\Gamma$, where the constant $C_N$ does not depend on $\gamma$. Here $|\gamma|$ is the length of $\gamma$ in the word metric on  $\Gamma$. Finally, the action of $\Gamma$ on $\mathcal{A}$ is required to be tempered: for any $m$ there exists $k$ and a real polynomial $P(z)$ such that $\|\gamma(a)\|_m\le P(|\gamma|)\|a\|_k$ for all $a$ and $\gamma$. 
The product in $\mathcal{A}\rtimes \Gamma$ is defined by the formula: 
$$
\{f_1(\gamma)\}\cdot \{f_2(\gamma)\}=\left\{\sum_{\gamma_1\gamma_2=\gamma}f_1(\gamma_1)\gamma_1(f_2(\gamma_2))\right\}.
$$

\begin{definition}
{\em The symbol of a $\Gamma$-operator} \eqref{eq-op1} is the pair $\sigma(\mathcal{D})=(\sigma_M(\mathcal{D}),\sigma_X(\mathcal{D}))$, consisting of the interior and the boundary symbols
\begin{equation}\label{s-bdm1}
 \sigma_M(\mathcal{D})=\{\sigma(A_\gamma)\}_{\gamma\in\Gamma}\in C^\infty(T^*_0M)\rtimes \Gamma, \quad \sigma_X(\mathcal{D})=\{\sigma_X(\mathcal{D}_\gamma)\}_{\gamma\in\Gamma}\in \Sigma_X\rtimes \Gamma.
\end{equation}
\end{definition}

\begin{definition}
A $\Gamma$-operator $\mathcal{D}$ is {\em elliptic}, if its interior and boundary symbols  are invertible in the corresponding crossed products \eqref{s-bdm1}.
\end{definition}
\begin{theorem}
An elliptic  $G$-operator has Fredholm property.
\end{theorem} 
The proof is standard (e.g., see~\cite{AnLe2}). More precisely, if $\mathcal{D}$ is elliptic, then its symbol is invertible. We denote the inverse symbol by 
$$
\left(\{b_{\gamma,M}\}_{\gamma\in\Gamma},\{b_{\gamma,X}\}_{\gamma\in\Gamma}\right)\in \Bigl(C^\infty(T^*_0M) \oplus \Sigma_X\Bigr)\rtimes \Gamma.
$$
Then, given $\gamma\in\Gamma,$ we choose a Boutet de Monvel operator $B_\gamma$ with the symbol $(b_{\gamma,M},b_{\gamma,X})$.
Finally, we define the operator  
$$
\mathcal{B}=\sum_{\gamma\in\Gamma} B_\gamma T_\gamma.
$$
A direct computation shows that $\mathcal{B}$ is a regularizer for $\mathcal{D}$ modulo compact operators.

The aim of this paper is to obtain an index formula for elliptic $\Gamma$-operators. To state our cohomological index 
formula, we note that the boundary of $T^*M$ is naturally fibered over $T^*X$.  In the next section, we describe de Rham complex on manifolds with fibered boundary and then define the components of the index formula in the cohomology of this complex.

\section{Two de Rham complexes for a manifold with fibered boundary}

In this section, let   $M$ be a Riemannian manifold with boundary $\partial M$ and suppose that the boundary is the total space of a locally-trivial fiber bundle $\pi:\partial M\to X$ with fiber $F$. Then the pair $(M,\pi)$ is called a {\em manifold with fibered boundary}. 

The embedding $i:\partial M\to M$ induces the restriction mapping $i^*:\Omega^*(M)\to \Omega^*(\partial M)$ of compactly supported differential forms to the submanifold. The projection $\pi$ defines the pull-back mapping $\pi^*:\Omega^*(X)\to \Omega^*(\partial M)$ and the mapping  
$$
\pi_*:\Omega^*(\partial M)\longrightarrow \Omega^{*-\nu}(X),\quad \nu=\dim F,
$$
which takes the forms to their integrals along the fibers of   $\pi$, e.g. see \cite{BGV1}. Here we suppose that  we are given a continuous family of orientations on the fibers of $\pi$.

Consider the morphism
$$
(\Omega^*(M),d) \stackrel{\pi_*i^*}\longrightarrow (\Omega^*(X),d)
$$
of de Rham complexes of   $M$ and $X$ and denote the cone of this homomorphism by $(\Omega^*(M,\pi),\partial)$, where
\begin{equation}\label{eq-c1}
 \Omega^j(M,\pi)=\Omega^j(M)\oplus\Omega^{j-\nu}(X), \quad
 \partial=\left(
            \begin{array}{cc}
             d & 0\\
             -\pi_*i^*& (-1)^\nu d
            \end{array}\right).
\end{equation}
The cohomology of the complex $(\Omega^*(M,\pi),\partial)$ is denoted by $H^* (M,\pi)$. 

Consider also the complex $(\widetilde{\Omega}^*(M,\pi),\widetilde\partial)$:
\begin{equation}\label{eq-c2}
 \widetilde{\Omega}^j(M,\pi)=\{(\omega,\omega_X)\in\Omega^j(M)\oplus\Omega^j(X)\;|\; i^*\omega=\pi^*\omega_X\}, \quad
 \widetilde{\partial}=\left(
            \begin{array}{cc}
             d & 0\\
             0& d
            \end{array}\right).  
\end{equation}
Denote the cohomology of $(\widetilde{\Omega}^*(M,\pi),\widetilde{\partial})$ by $\widetilde{H}^*(M,\pi)$.

Componentwise products of differential forms  give us the product
$$
\wedge: \Omega^j(M,\pi)\times \widetilde{\Omega}^k(M,\pi)\longrightarrow \Omega^{j+k}(M,\pi).
$$
The following Leibniz rule
$$
\partial (a\wedge b)=\partial a\wedge b+(-1)^{j}a\wedge \widetilde\partial b,\quad a\in \Omega^j(M,\pi), b\in \widetilde{\Omega}^k(M,\pi) 
$$
implies that the product $\wedge$ defines a product in cohomology
$$
 \wedge:H^j(M,\pi)\times \widetilde{H}^k(M,\pi)\longrightarrow H^{j+k}(M,\pi).
$$

Finally, suppose that   $M$ and $X$ are oriented and the orientation in the fibers of  $\pi$ is compatible with the orientations of the total space and the base. Then we have a well defined mapping   
$$
\begin{array}{rcc}
  \displaystyle\langle \cdot,[M,\pi]\rangle:H^* (M,\pi) & \longrightarrow & \mathbb{C}\vspace{2mm}\\
   (\omega,\omega_X) & \longrightarrow &  \displaystyle\int_M \omega+\int_X \omega_X.
\end{array}
$$

\section{Index formula}

\paragraph{Noncommutative differential forms. Regularized traces.}

Denote by $\widetilde\Sigma_X$ the algebra of boundary symbols, which are defined on $T^*X$ and twisted homogeneous at infinity (i.e., for large $|\xi'|$).
 
Consider actions of $\Gamma$ on the algebras $\Omega(T^*M),\widetilde\Sigma_X\otimes_{C^\infty(X)}\Omega(T^*X)  $ of compactly supported differential forms and the corresponding smooth crossed products
$$
\Omega(T^*M)\rtimes \Gamma, \quad \left(\widetilde\Sigma_X\otimes_{C^\infty(X)}\Omega(T^*X)\right)\rtimes \Gamma.
$$
These products are differential graded algebras.

Given $\gamma\in \Gamma$, we define mappings (cf. \cite{NaSaSt17,SaSt26})
\begin{equation}\label{eq-tr1}
 \tau^{\gamma}:\Omega(T^*M)\rtimes \Gamma \longrightarrow \Omega(T^*M^{\gamma}),
\end{equation}
\begin{equation}\label{eq-tr2}
 \tau^\gamma_{X}:\left(\widetilde\Sigma_X\otimes_{C^\infty(X)}\Omega(T^*X)\right)\rtimes \Gamma \longrightarrow \Omega(T^*X^{\gamma}).
\end{equation}
To define these mappings, we introduce some notation. Denote by  $\overline{\Gamma}$ the closure of $\Gamma$ in the compact Lie group of isometries of $M$. This closure is a compact Lie group. Let $C^\gamma\subset\overline{\Gamma}$ be the centralizer \footnote{Recall that the centralizer of  $\gamma$ is the subgroup of elements commuting with $\gamma$.}
of $\gamma$. The centralizer is a closed Lie subgroup in   $\overline{\Gamma}$. The elements of the centralizer are denoted by  $h$, while the induced Haar measure on the centralizer is denoted by $dh$. Below, given $\gamma'\in \langle
\gamma\rangle$ we fix an arbitrary element $z=z(\gamma,\gamma')$ which conjugates $\gamma$ and
$\gamma'=z\gamma z^{-1}$. Any such element defines a diffeomorphism  $\partial z:T^*M^\gamma\to T^*M^{\gamma'}$ of the corresponding fixed point sets.

We define the functional  \eqref{eq-tr1}  as
\begin{equation}\label{eq-sled-nash1}
    \tau^\gamma(\omega)=
      \sum_{\gamma'\in \langle \gamma\rangle}\;\;\;
        \int_{C^\gamma}
           \Bigl.h^*\bigl(
              {z}^*\omega(\gamma')
             \bigr)\Bigr|_{T^*M^\gamma}
           dh,\quad \text{where }\omega\in \Omega(T^*M)\rtimes \Gamma,
\end{equation}
while the functional  \eqref{eq-tr2} as
\begin{equation}\label{eq-sled-nash2}
    \tau^\gamma_X(\omega_X)=
      \sum_{\gamma'\in \langle \gamma\rangle}\;\;\;
        \int_{C^\gamma} \tr_X
           \Bigl.h^*\bigl(
              {z}^*\omega_X(\gamma')
             \bigr)\Bigr|_{T^*X^\gamma}
           dh,\quad \text{where }\omega_X\in \left(\widetilde\Sigma_X\otimes_{C^\infty(X)}\Omega(T^*X)\right)\rtimes \Gamma.
\end{equation}
Here 
$$
\tr_X \left(
\sum_I \omega_I(z)dz^I\right)= \sum_I \tr'(\omega_I(z))dz^I,
$$
where the regularized trace $\tr':\widetilde\Sigma_X \to C^\infty(T^*X)$ was defined earlier.

One can show that the expressions   \eqref{eq-sled-nash1} and \eqref{eq-sled-nash2} do not depend on the choice of $z$.

Clearly, these functionals enjoy the properties
\begin{equation}\label{eq-trcl1}
\tau^\gamma(\omega_1\wedge\omega_2)=(-1)^{\deg\omega_1\deg\omega_2}\tau^\gamma(\omega_2\wedge\omega_1),\quad \omega_1,\omega_2\in \Omega(T^*M)\rtimes\Gamma
\end{equation}
\begin{equation}\label{eq-trcl2}
d\tau^\gamma_X(\omega )= \tau^\gamma_X(d\omega ),\quad \omega \in \left(\widetilde\Sigma_X\otimes_{C^\infty(X)}\Omega(T^*X)\right)\rtimes \Gamma,
\end{equation}
i.e., $\tau^\gamma$ is a graded trace, while $\tau_X^\gamma$ is closed.

\paragraph{Chern character of elliptic symbols.}

Let $\mathcal{D}$ be an $N\times N$ matrix elliptic $\Gamma$-operator. Then its interior and boundary symbols are invertible elements in the corresponding crossed products  and we denote the inverse symbols by
\begin{equation}\label{s-bdm2}
 \sigma_M(\mathcal{D})^{-1}\in C^\infty(T^*_0M,\Mat_N)\rtimes \Gamma, \quad \sigma_X(\mathcal{D})^{-1}\in (\Sigma_X\otimes\Mat_N) \rtimes \Gamma.
\end{equation}
Let us extend $\sigma_M(\mathcal{D})^{\pm 1}$ to $T^*M$ up to a smooth symbol satisfying transmission property and homogeneous at infinity, and extend 
$\sigma_X(\mathcal{D})^{\pm 1}$ to $T^*X$ as a smooth symbol twisted homogeneous at infinity. Denote these extensions by
$$
a,r\in C^\infty(T^*M)\rtimes \Gamma, \qquad a_X,r_X\in  \left(\widetilde\Sigma_X\otimes_{C^\infty(X)}\Omega(T^*X)\right)\rtimes \Gamma.
$$
We shall suppose that these extensions are compatible, i.e., the symbol of the boundary symbol is equal to the restriction of the interior symbol to the boundary.   

Let us define noncommutative connections in the trivial rank  $N$ bundles over $T^*M$ and $T^*X$  as
$$
\nabla_M= d+rda\wedge, \quad \nabla_X= d+r_Xda_X\wedge.
$$
Their curvature forms are equal to 
$$
\Omega_M=\nabla_M^2=dr\wedge da+(rda)^2,\quad \Omega_X=\nabla_X^2=dr_X\wedge da_X+(r_Xda_X)^2.
$$
Let us define compactly supported differential forms  
\begin{equation}\label{eq-forms1}
 \ch_{T^*M}^\gamma\sigma(\mathcal{D}) \in \Omega^{ev}(T^*M^\gamma),\quad \ch^\gamma_{T^*X}\sigma_X(\mathcal{D})\in \Omega^{ev}(T^*X^\gamma)
\end{equation}
on the cotangent bundles of the fixed point submanifolds as  
$$ 
\ch^\gamma_{T^*M}\sigma(\mathcal{D})=\tau^{\gamma}\left(\exp\left(-\frac{\Omega_M}{2\pi i}\right)(1_N-ra)\right)- 
\tau^{\gamma}\left(1_N-a\exp\left(-\frac{\Omega_M}{2\pi i}\right)r\right)
$$
$$ 
\ch^\gamma_{T^*X}\sigma(\mathcal{D})=\tau^\gamma_X\left(\exp\left(-\frac{\Omega_X}{2\pi i}\right)(1_N-r_Xa_X)\right)- 
\tau^\gamma_X\left(1_N-a_X\exp\left(-\frac{\Omega_X}{2\pi i}\right)r_X\right).
$$ 
Since    $\tau^\gamma$ is a graded trace, the first form can be written in a simpler way  
$$
\ch^\gamma_{T^*M}\sigma(\mathcal{D})=\tau^{\gamma}\left(\exp\left(-\frac{\Omega_M}{2\pi i}\right) \right)-N.
$$

The boundary $\partial(T^*M^\gamma)\simeq T^*X^\gamma\times\mathbb{R}$ is fibered over $T^*X^\gamma$ with fiber $\mathbb{R}$.
Denote the corresponding projection by $\pi^\gamma:\partial(T^*M^\gamma)\to T^*X^\gamma$ and the embedding
$\partial (T^*M^\gamma)\subset T^*M^\gamma$ by $i_\gamma$.
\begin{proposition}
Given $\gamma\in\Gamma$, the pair $(\ch_{T^*M}^\gamma\sigma(\mathcal{D}),\ch_{T^*X}^\gamma\sigma_X(\mathcal{D}))$ enjoys the properties
\begin{equation}\label{eq-fund1}
d\left(\ch_{T^*M}^\gamma\sigma(\mathcal{D})\right)=0,\quad d\left(\ch_{T^*X}^\gamma\sigma(\mathcal{D})\right)=\pi^\gamma_*i_\gamma^*\left(\ch_{T^*M}^\gamma\sigma(\mathcal{D})\right),
\end{equation}
i.e., it is closed in the complex $(\Omega^*(T^*M^\gamma,\pi^\gamma),\partial)$, see \eqref{eq-c1}, and its cohomology class, denoted by
$$
 \ch^\gamma\sigma(\mathcal{D}) \in H^{ev} (T^*M^\gamma,\pi^\gamma),
$$
does not depend on the choice of $a,r,a_X,r_X$ and does not change under homotopies of elliptic symbols.  
\end{proposition}
 
\paragraph{Index formula}

To state the index formula, we define the necessary equivariant characteristic classes.
First, we define the Todd form  on $M^\gamma$:
$$
\Td(T^*M^\gamma\otimes\mathbb{C})=\det\left(
\frac{-\Omega^\gamma/2\pi i}{1-\exp(\Omega^\gamma/2\pi i)}\right),
$$
where $\Omega^\gamma$ is the curvature form of the Levi-Civita connection on $M^\gamma$.
One similarly defines the Todd form  $\Td(T^*X^\gamma\otimes\mathbb{C})$ on $X^\gamma$. The pair of these forms is closed in the complex  $(\widetilde{\Omega}^*(M^\gamma,\pi^\gamma),\widetilde\partial)$, see \eqref{eq-c2}, its cohomology class is denoted by
$$
\Td^\gamma(T^*M\otimes\mathbb{C}) \in \widetilde{H}^{ev}(M^\gamma,\pi^\gamma).
$$
Second, let $N^\gamma$ be the normal bundle of $M^\gamma\subset M$. Then we have a natural action of $\gamma$ on $N^\gamma$ and the following differential form on $M^\gamma$ is defined:
$$
 \ch^\gamma\Lambda( {N}^\gamma\otimes \mathbb{C})=\tr_{\Lambda^{ev}(N^\gamma)}\left(\gamma\exp(-\Omega/2\pi i)\right)- \tr_{\Lambda^{odd}(N^\gamma)}\left(\gamma\exp(-\Omega/2\pi i)\right),
$$
where $\Omega$ is the curvature form of the exterior bundle $\Lambda(N^\gamma)$, $\gamma$ is considered  as an endomorphism of the subbundles $\Lambda^{ev/odd}(N^\gamma)$  of even/odd forms and 
$ 
\tr_{\Lambda^{ev/odd}(N^\gamma)} 
$ 
is the fiberwise trace functional of endomorphisms of $\Lambda^{ev/odd}(N^\gamma)$. 
Similarly, let $N^\gamma_X$ be the normal bundle of $X^\gamma\subset X$. then one defines the form $\ch^\gamma\Lambda( {N}_X^\gamma\otimes \mathbb{C})$ on $X^\gamma$.
The pair $(\ch^\gamma\Lambda( {N}^\gamma\otimes \mathbb{C}),\ch^\gamma\Lambda( {N}_X^\gamma\otimes \mathbb{C}))$ is closed in the complex \eqref{eq-c2}. We denote its cohomology class by  
$$
 \ch^\gamma\Lambda(\mathcal{N}^\gamma\otimes \mathbb{C})\in  \widetilde{H}^{ev}(M^\gamma,\pi^\gamma).
$$
The latter class is invertible, since its zero degree component is a nonzero complex number (see~\cite{AtSi1} or \cite{NaSaSt17} for a proof).

\begin{theorem}The following index formula holds:
\begin{equation}\label{eq-indf1}
 \ind \mathcal{D}=\sum_{\langle \gamma\rangle\subset \Gamma}\langle \ch^\gamma\sigma(\mathcal{D})\wedge \Td^\gamma(T^*M\otimes\mathbb{C})\wedge
 \ch^\gamma\Lambda(\mathcal{N}^\gamma\otimes \mathbb{C})^{-1},[T^*M^\gamma,\pi^\gamma]\rangle
\end{equation}
\end{theorem}


\end{document}